\documentclass[12pt]{article}
\usepackage{amsmath,eufrak}
\newcommand{\cl}{C\kern -0.2em \ell}
\newcommand{\BH}{{\rm\hskip 0.1pt %
                         I\hskip -2.15pt H}}
\newcommand{\K}{\mbox{\bf K}}
\newcommand{\R}{\mbox{\bf R}}
\newcommand{\C}{\mbox{\bf C}}
\newcommand{\Z}{\mbox{\bf Z}}
\newcommand{\e}{\mbox{\bf e}}
\newcommand{\M}{\mbox{\rm M}}
\newcommand{\der}{\mbox{\rm d}}
\newcommand{\re}{\mbox{\rm Re}\,}
\newcommand{\im}{\mbox{\rm Im}\,}

\newcommand{\ov}{\overline}
\DeclareMathOperator{\Spin}{\bf Spin}
\begin{document}
\begin{center}
{\LARGE Generalized Weierstrass representation for\protect\newline 
surfaces in terms of Dirac-Hestenes\\[0.65cm] spinor field}\\[1.5cm]
{\large Vadim V. Varlamov}
\\[0.3cm]
{\small\it Computer Division, Siberia State University of Industry,}\\
{\small\it Novokuznetsk 654007, Russia}
\vspace{1.5cm}
\begin{abstract}
A representation of generalized Weierstrass formulae 
for an immersion of generic
surfaces into a 4-dimensional complex space in terms of spinors treated 
as minimal
left ideals of Clifford algebras is proposed. 
The relation between integrable
deformations of surfaces via mVN-hierarchy and integrable deformations of
spinor fields on the surface is also discussed.
\end{abstract}
\end{center}
\hfill

\underline{Mathematics Subject Classification (1991)}: 53A05, 53A10, 15A66.

\underline{Key words}: Weierstrass representation, Clifford algebras, spinors.
\setcounter{page}{0}
\thispagestyle{empty}
\newpage
\section{Introduction}
The theory of integrable deformations and immersions of surfaces due its a
close relationship with the theory of integrable systems at present time
is a rapid developing area of mathematical physics. One of the most
powerful methods in this area is a Weierstrass representation for minimal
surfaces \cite{1}, the generalization of which onto a case of generic
surfaces was proposed by Konopelchenko in 1993 \cite{2,3} served as
a basis for the following investigations. So, the generalized Weierstrass
formulae for conformal immersion of surfaces into 3-dimensional Euclidean
space are used for the study of the basic quantities related to 2D gravity,
such as Polyakov extrinsic action, Nambu-Goto action, geometric action
and Euler characteristic \cite{4}. This method is also intensively used
for the study of constant mean curvature surfaces, Willmore surfaces,
surfaces of revolution and in many other problems related with
differential geometry \cite{5}-\cite{14}. 
A further generalization of Weierstrass
representation onto a case of multidimensional Riemann spaces, in particular
onto a case of 4-dimensional space with signature $(+,+,+,-)$ (Minkowski
space-time) has been proposed in the recent paper \cite{15}.

In the present paper we consider a relation between a Weierstrass representation
in a 4-dimensional complex space $\C^4$ and a Dirac-Hestenes spinor field which
is defined in Minkowski space-time $\R^{1,3}$. Dirac-Hestenes spinors 
were originally 
introduced
in \cite{16,17} for the formulation of a Dirac theory of electron with the
usage of the
space-time algebra $\cl_{1,3}$ \cite{18} in $\R^{1,3}$ (see also \cite{19}).
On the other hand, there is a very graceful formulation \cite{20}-\cite{23}
of the Dirac-Hestenes theory in terms of modern interpretation of spinors
as minimal left ideals of Clifford algebras \cite{24,25}, a brief review
of which we give in Section 2. In Section 3 after a short historical
introduction, generalized Weierstrass formulae in $\C^4$ are rewritten
in a spinor representation type form (matrix representation of a biquaternion
algebra $\C_2\cong\M_2(\C)$) and are identified with the Dirac-Hestenes
spinors, the matrix representation of which is also isomorphic to
$\M_2(\C)$. It allows to use a well-known relation between Dirac-Hestenes  
and Dirac spinors \cite{23,26} (see also \cite{27}) to establish a relation
between Weierstrass-Konopelchenko coordinates for surfaces immersed into
$\C^4$ and Dirac spinors. Integrable deformations of surfaces defined
by a modified Veselov-Novikov equation and their relation with  
integrable deformations of Dirac field on surface
are considered at the end of the Section 3. 
\section{Spinors as minimal left ideals of Clifford algebras}
Let us consider a Clifford algebra $\cl_{p,q}(V,Q)$ over a field $\K$ of
characteristic 0 $(\K=\R,\,\K=\Omega=\R\oplus\R,\,\K=\C)$, where $V$ is a
vector space endowed with a nondegenerate quadratic form
$$Q=x^{2}_{1}+\ldots+x^{2}_{p}-\ldots-x^{2}_{p+q}.$$
The algebra $\cl_{p,q}$ is naturally $\Z_{2}$-graded. Let $\cl^{+}_{p,q}$
(resp. $\cl^-_{p,q}$) 
be a set consisting of all even (resp. odd)
elements of algebra $\cl_{p,q}$.
The set $\cl^{+}_{p,q}$ is a
subalgebra of $\cl_{p,q}$. It is obvious that $\cl_{p,q}=\cl^{+}_{p,q}
\oplus\cl^{-}_{p,q}$.

When $n$ is odd, a volume element $\omega=\e_{12\ldots p+q}$ commutes with all
elements of algebra $\cl_{p,q}$ and therefore belongs to a center of
$\cl_{p,q}$. Thus, in the case of $n$ is odd we have for a center
\begin{equation}\label{e1}
\Z_{p,q}=\begin{cases}
\R\oplus i\R & \text{if $\omega^{2}=-1$};\\
\R\oplus e\R & \text{if $\omega^{2}=+1$},
\end{cases}\end{equation}
where $e$ is a double unit. In the case of $n$ is even the center of
$\cl_{p,q}$ consists the unit of algebra.

Let $\R_{p,q}=\cl_{p,q}(\R^{p,q},Q)$ be a real Clifford algebra ($V=\R^{p,q}$
is a real space). Analogously, in the case of a complex space we have 
$\C_{p,q}=\cl_{p,q}(\C^{p,q},Q)$. 
 Moreover, it is obvious that $\C_{p,q}\cong\C_{n}$, where
$n=p+q$. Further, let us consider the following most important in physics
Clifford algebras and their isomorhisms to matrix algebras:
\begin{eqnarray}
\text{quaternions}&&\R_{0,2}=\BH \nonumber\\
\text{biquaternions}&&\C_{2}=\R_{3,0}\cong\M_{2}(\C)\nonumber\\
\text{space-time algebra}&&\R_{1,3}\cong\M_{2}(\BH)\nonumber\\
\text{Dirac algebra}&&\C_{4}=\R_{4,1}\cong\M_{4}(\C)\cong\M_{2}(\C_{2})
\nonumber
\end{eqnarray}
The identity $\C_{2}=\R_{3,0}$ for a biquaternion algebra known in physics as
a Pauli algebra is immediately obtained from definition of the center of the
algebra $\cl_{p,q}$ (\ref{e1}). Namely, for $\R_{3,0}$ we have a volume element
$\omega=\e_{123}\in\Z_{3,0}=\R\oplus i\R$, since $\omega^{2}=-1$. The
identity $\C_{4}=\R_{4,1}$ is analogously proved. The isomorphism
$\R_{4,1}\cong\M_{2}(\C_{2})$ is a consequence of an algebraic modulo 2
periodicity of complex Clifford algebras: $\C_{4}\cong\C_{2}\otimes\C_{2}
\cong\C_{2}\otimes\M_{2}(\C)\cong\M_{2}(\C_{2})$ \cite{28,29,30}.  

The left (resp. right) ideal of algebra $\cl_{p,q}$ is defined by the
expression $\cl_{p,q}e$ (resp. $e\cl_{p,q}$), where $e$ is an idempotent
satisfying the condition $e^{2}=e$. Analogously, a minimal left (resp. right)
ideal is a set of type $I_{p,q}=\cl_{p,q}e_{pq}$ (resp. $e_{pq}\cl_{p,q}$),
where $e_{pq}$ is a primitive idempotent, i.e., $e_{pq}^{2}=e_{pq}$ and
$e_{pq}$ cannot be represented as a sum of two orthogonal idempotents, i.e.,
$e_{pq}\neq f_{pq}+g_{pq}$, where $f_{pq}g_{pq}=g_{pq}f_{pq}=0,\;f^{2}_{pq}=
f_{pq},\;g^{2}_{pq}=g_{pq}$. In the general case a primitive idempotent
has a form \cite{20}
\begin{equation}\label{e2}
e_{pq}=\frac{1}{2}(1+\e_{\alpha_{1}})\frac{1}{2}(1+\e_{\alpha_{2}})\ldots
\frac{1}{2}(1+\e_{\alpha_{k}}),
\end{equation}
where $\e_{\alpha_{1}},\ldots,\e_{\alpha_{k}}$ are commuting elements of the
canonical basis of $\cl_{p,q}$ such that $(\e_{\alpha_{i}})^{2}=1,\;
(i=1,2,\ldots,k)$. The values of $k$ are defined by a formula
\begin{equation}\label{e3}
k=q-r_{q-p},
\end{equation}
where $r_{i}$ are the Radon-Hurwitz numbers, values of which form
a cycle of the period 8 :
\begin{equation}\label{e4}
r_{i+8}=r_{i}+4.
\end{equation}
The values of all $r_{i}$ are
\begin{center}
\begin{tabular}{lcccccccc}
$i$   & 0 & 1 & 2 & 3 & 4 & 5 & 6 & 7 \\ \hline
$r_i$ & 0 & 1 & 2 & 2 & 3 & 3 & 3 & 3 \\
\end{tabular}
\end{center}
For example, let consider a minimal left ideal of the space-time algebra 
$\R_{1,3}$.
The Radon-Hurwitz number for algebra $\R_{1,3}$ is equal to $r_{q-p}=r_2=2$,
and therefore from (\ref{e3}) we have $k=1$. The primitive idempotent of
$\R_{1,3}$ has a form
$$e_{13}=\frac{1}{2}(1+\e_{0}),$$
or $e_{13}=\frac{1}{2}(1+\Gamma_{0})$, where $\Gamma_{0}$ is a matrix
representation of the unit $\e_{0}\in\R_{1,3}$. Thus, a minimal left ideal
of $\R_{1,3}$ is defined by the following expression
\begin{equation}\label{e5}
I_{1,3}=\R_{1,3}\frac{1}{2}(1+\Gamma_{0}).
\end{equation}
Analogously, for the Dirac algebra $\R_{4,1}$ on using the 
recurrence formula
(\ref{e4}) we obtain $k=1-r_{-3}=1-(r_5-4)=2$, and
a primitive idempotent of $\R_{4,1}$ may be defined as follows
\begin{equation}\label{e6}
e_{41}=\frac{1}{2}(1+\Gamma_{0})\frac{1}{2}(1+i\Gamma_{12}),
\end{equation}
where $\Gamma_{12}=\Gamma_{1}\Gamma_{2}$ and $\Gamma_{i}\;(i=0,1,2,3)$ are
matrix representations of the units of $\R_{4,1}=\C_4$:
$$\Gamma_{0}=
{\renewcommand{\arraystretch}{1}
\begin{pmatrix}
1 & 0 & 0 & 0 \\
0 & 1 & 0 & 0 \\
0 & 0 &-1 & 0 \\
0 & 0 & 0 &-1
\end{pmatrix},\quad\Gamma_{1}=
\begin{pmatrix}
0 & 0 & 0 & 1 \\
0 & 0 & 1 & 0 \\
0 &-1 & 0 & 0 \\
-1 & 0 & 0 & 0
\end{pmatrix}},$$$$\Gamma_{2}=
{\renewcommand{\arraystretch}{1}
\begin{pmatrix}
0 & 0 & 0 &-i \\
0 & 0 & i & 0 \\
0 & i & 0 & 0 \\
-i & 0 & 0 & 0
\end{pmatrix},\quad\Gamma_{3}=
\begin{pmatrix}
0 & 0 & 1 & 0 \\
0 & 0 & 0 &-1 \\
-1 & 0 & 0 & 0 \\
0 & 1 & 0 & 0
\end{pmatrix}}.$$

Further, for a minimal left ideal of Dirac algebra $I_{4,1}=\R_{4,1}
\frac{1}{2}(1+\Gamma_{0})\frac{1}{2}(1+i\Gamma_{12})$ using the isomorphisms
$\R_{4,1}=\C_{4}=\C\otimes\R_{1,3}\cong\M_{2}(\C_{2}),\;\R^{+}_{4,1}\cong
\R_{1,3}\cong\M_{2}(\BH)$ and also an identity $\R_{1,3}e_{13}=\R^{+}_{1,3}
e_{13}$ \cite{22,23} we have the following expression \cite{27}:
\begin{multline}\label{e7}
I_{4,1}=\R_{4,1}e_{41}=(\C\otimes\R_{1,3})e_{41}\cong\R^{+}_{4,1}e_{41}
\cong\R_{1,3}e_{41}=\\
R_{1,3}e_{13}\frac{1}{2}(1+i\Gamma_{12})=\R^{+}_{1,3}e_{13}\frac{1}{2}
(1+i\Gamma_{12}).
\end{multline}

Let $\Phi\in\R_{4,1}\cong\M_{4}(\C)$ be a Dirac spinor and $\phi\in\R^{+}_{1,3}
\cong\R_{3,0}=\C_{2}$ be a Dirac-Hestenes spinor. Then from (\ref{e7})
the relation immediately follows  between spinors $\Phi$ and $\phi$:
\begin{equation}\label{e8}
\Phi=\phi\frac{1}{2}(1+\Gamma_{0})\frac{1}{2}(1+i\Gamma_{12}).
\end{equation}
Since $\phi\in\R^{+}_{1,3}\cong\R_{3,0}$, then the Dirac-Hestenes spinor
can be represented by a biquaternion number
\begin{equation}\label{e9}
\phi=a^{0}+a^{01}\Gamma_{01}+a^{02}\Gamma_{02}+a^{03}\Gamma_{03}
+a^{12}\Gamma_{12}+a^{13}\Gamma_{13}+a^{23}\Gamma_{23}+a^{0123}\Gamma_{0123}.
\end{equation}
Or in the matrix representation
\begin{equation}\label{e10}
\phi=
{\renewcommand{\arraystretch}{1}
\begin{pmatrix}
\phi_{1} & -\phi^{\ast}_{2} & \phi_{3} & \phi^{\ast}_{4} \\
\phi_{2} & \phi^{\ast}_{1} & \phi_{4} & -\phi^{\ast}_{3} \\
\phi_{3} & \phi_{4}^{\ast} & \phi_{1} & -\phi^{\ast}_{2} \\
\phi_{4} & -\phi^{\ast}_{3} & \phi_{2} & \phi^{\ast}_{1}
\end{pmatrix}},\quad
\phi_{i}\in\C,
\end{equation}
where
\[
\phi_1=a^0-ia^{12},\quad
\phi_2=a^{31}-ia^{23},\quad
\phi_3=a^{03}-ia^{0123},\quad
\phi_4=a^{01}+ia^{02}.
\]
Finally, from (\ref{e8}) it follows that for the Dirac spinor $\Phi$ and also
a space-time spinor $Z=\phi\frac{1}{2}(1+\Gamma_{0})$ we have expressions
$$\Phi=
{\renewcommand{\arraystretch}{1}
\begin{pmatrix}
\phi_{1} & 0 & 0 & 0 \\
\phi_{2} & 0 & 0 & 0 \\
\phi_{3} & 0 & 0 & 0 \\
\phi_{4} & 0 & 0 & 0
\end{pmatrix},\quad
Z=
\begin{pmatrix}
\phi_{1} & -\phi_{2}^{\ast} & 0 & 0 \\
\phi_{2} & \phi^{\ast}_{1} & 0 & 0 \\
\phi_{3} & \phi^{\ast}_{4} & 0 & 0 \\
\phi_{4} & -\phi^{\ast}_{3} & 0 & 0
\end{pmatrix}},$$
which are minimal left ideals of algebras $\R_{4,1}$ and $\R_{1,3}$,
respectively.

The Dirac spinor $\Phi$ may be considered as a vector in the 4-dimensional
complex space $\C^4$ associated with the algebra $\C_4$. However, from a
physical point of view it is more natural to consider the spinor $\Phi$ in
space-time $\R^{1,3}$. In connection with this, let us introduce
(follows \cite{21,22,23,26}) a more rigorous definition of spinor as a
minimal left ideal of algebra $\cl_{p,q}$.

Let $\mathfrak{B}_\Sigma=\left\{\Sigma_0,\overset{.}{\Sigma},\overset{..}
{\Sigma},\ldots\right\}$ be a set of all ordered orthonormal bases for
$\R^{p,q}$. Any two bases $\Sigma_0,\,\overset{.}{\Sigma}\in
\mathfrak{B}_\Sigma$ are related by the element of the group $\Spin_+(p,q)$:
\[
\overset{.}{\Sigma}=u\Sigma_0u^{-1},\quad u\in\Spin_+(p,q).
\]
Analogously, for the primitive idempotents defined in the basis $\Sigma\in
\mathfrak{B}_\Sigma$ and denoted as $e_{\Sigma_0},e_{\overset{.}{\Sigma}},
\ldots$, we have $e_{\overset{.}{\Sigma}}=ue_{\Sigma_0}u^{-1},\;u\in
\Spin_+(p,q)$. Then the ideals $I_{\Sigma_0}, I_{\overset{.}{\Sigma}},
I_{\overset{..}{\Sigma}},\ldots$ are geometrically equivalent if and only
if
\[
I_{\overset{.}{\Sigma}}=uI_{\Sigma_0}u^{-1},\quad u\in\Spin_+(p,q),
\]
or, since $uI_{\Sigma_0}=I_{\Sigma_0}$:
\[
I_{\overset{.}{\Sigma}}=I_{\Sigma_0}u^{-1}.
\]
Therefore, an algebraic spinor for $\R^{p,q}$ is an equivalence class of the
quotent set $\left\{I_{\Sigma}\right\}/\R$, where $\left\{I_\Sigma\right\}$
is a set of all geometrically equivalent ideals, and $\Phi_{\Sigma_0}\in
I_{\Sigma_0}$ and $\Phi_{\overset{.}{\Sigma}}\in I_{\overset{.}{\Sigma}}$
are equivalent, $\Phi_{\overset{.}{\Sigma}}\cong\Phi_{\Sigma_0}\pmod{\R}$
if and only if $\Phi_{\overset{.}{\Sigma}}=\Phi_{\Sigma_0}u^{-1},\;
u\in\Spin_+(p,q)$.   
\section{Weierstrass representation for surfaces in space $\C^4$}
Historically, the Weierstrass representation \cite{1} appeared in the result 
of the
following variational problem: among the surfaces restricted by the some
curve for finding such a surface, the area of which is minimal, i.e., it is
necessary to find a minimum of the functional
$$S=\int\int\sqrt{1+p^{2}+q^{2}}dxdy,$$
where $p=dz/dx,\;
q=dz/dy,\;z=f(x,y)$ is an
equation of the surface. The Euler equation for this problem has a form
$$\frac{\partial}{\partial x}\left(\frac{p}{\sqrt{1+p^2+q^2}}\right)+
\frac{\partial}{\partial y}\left(\frac{q}{\sqrt{1+p^2+q^2}}\right)=0.$$
This equation expresses a main geometrical property of such a surface :
in each point the meant curvature is equal to zero. The surface 
which possesses 
such a property is called {\it a minimal surface}. If we compare a region
$\mathfrak{M}$ of the surface with a region $\mathfrak{E}$ of the flat surface
so that the point on $\mathfrak{M}$ with the coordinates $(X^1,\,X^2,\,X^3)$
corresponds to a point $w=u+iv$ of region $\mathfrak{E}$, then for the
minimal surface we have the equations
$$\frac{\partial^{2}X^1}{\partial u^2}+\frac{\partial^2X^1}{\partial v^2}=0,
\quad\frac{\partial^2X^2}{\partial u^2}+\frac{\partial^2X^2}{\partial v^2}=0,
\quad\frac{\partial^2X^3}{\partial u^2}+\frac{\partial^2X^3}{\partial v^2}=0,$$
solutions of which are of the form
$$X^1=\re f(w),\quad X^2=\re g(w),\quad X^3=\re h(w),$$
at
$$(f^{\prime}(w))^2+(g^{\prime}(w))^2+(h^{\prime}(w))^2=0.$$
The functions satisfying this equation are
$$f^{\prime}(w)=i(G^2+H^2),\quad g^{\prime}(w)=G^2-H^2,\quad h^{\prime}(w)=
2GH,$$
whence
\begin{eqnarray}
X^1&=&C^1+\re\int^{w}_{w_0}i(G^2+H^2)\der w,\nonumber \\
X^2&=&C^2+\re\int^{w}_{w_0}(G^2-H^2)\der w, \nonumber \\
X^3&=&C^3+2\re\int^{w}_{w_0}GH\der w.\label{e11} 
\end{eqnarray}
Here $G(w)$ and $H(w)$ are holomorphic functions defined in a circle or in 
all complex plane. After substitution of variables
$$s=\xi+i\eta=\frac{H(w)}{G(w)},\quad G^2\frac{\der w}{\der s}=F(s),$$
the equations (\ref{e11}) take the form
\begin{eqnarray}
\der X^1&=&\re\left[i(1+s^2)F(s)\der s\right],\nonumber\\
\der X^2&=&\re\left[(1-s^2)F(s)\der s\right],\nonumber\\
\der X^3&=&\re\left[2sF(s)\der s\right].\nonumber
\end{eqnarray}
Thus, for an every analytic function $F(s)$ we have a minimal surface.

Further, let us consider generalized Weierstrass representation for surfaces
immersed into 4-dimensional complex space $\C^4$, which, as known, is 
associated with the Dirac algebra $\C_4$. In this case generalized Weierstrass
formulae have a form
\begin{eqnarray}
X^1&=&\frac{i}{2}\int_\Gamma(\psi_1\psi_2d\ov{z}-\varphi_1\varphi_2dz),
\nonumber\\
X^2&=&\frac{1}{2}\int_\Gamma(\psi_1\psi_2d\ov{z}+\varphi_1\varphi_2dz),
\nonumber\\
X^3&=&\frac{1}{2}\int_\Gamma(\psi_1\varphi_2d\ov{z}-\varphi_1\psi_2dz),
\nonumber\\
X^4&=&\frac{i}{2}\int_\Gamma(\psi_1\varphi_2d\ov{z}+\varphi_1\psi_2dz),
\label{e16'}
\end{eqnarray}
\begin{equation}\label{e17''}
\begin{array}{ccc}
\psi_{\alpha z}&=&p\varphi_\alpha,\\
\varphi_{\alpha\ov{z}}&=&-p\psi_\alpha
\end{array}\quad,\alpha=1,2,
\end{equation}
where $\psi,\varphi$ and $p$ are complex-valued functions on variables
$z,\ov{z}\in\C$, $\Gamma$ is a contour in complex plane $\C$. 
We will interpret the functions $X^i(z,\ov{z})$ as the coordinates in
$\C^4$. It is easy to verify that components of an induced metric have
a form
\begin{eqnarray}
g_{zz}&=&\ov{g_{\bar{z}\bar{z}}}=\sum^4_{i=1}(X^i_z)^2=0,\nonumber\\
g_{z\bar{z}}&=&\sum^4_{i=1}(X^i_zX^i_{\bar{z}})=\psi_1\psi_2\varphi_1
\varphi_2.\nonumber
\end{eqnarray}
Therefore, the formulae
(\ref{e16'}), (\ref{e17''}) define a conformal immersion of the surface into
$\C^4$ with an induced metric
$$ds^2=\psi_1\psi_2\varphi_1\varphi_2dzd\ov{z}.$$

The formulae (\ref{e16'}) may be rewritten in the following form
\begin{eqnarray}
d(X^1+iX^2)&=&i\psi_1\psi_2\ov{z},\nonumber \\
d(X^1-iX^2)&=&-i\varphi_1\varphi_2dz,\nonumber \\
d(X^4+iX^3)&=&i\psi_1\varphi_2\ov{z},\nonumber\\
d(X^4-iX^3)&=&i\varphi_1\psi_2dz,\nonumber
\end{eqnarray}
or
\begin{equation}\label{e16}
d(X^4\sigma_0+X^1\sigma_1+X^2\sigma_2+X^3\sigma_3)=i
{\renewcommand{\arraystretch}{1}
\begin{pmatrix}
\varphi_1\psi_2dz & \psi_1\psi_2d\ov{z} \\
\varphi_1\varphi_2dz & \psi_1\varphi_2d\ov{z}
\end{pmatrix}},
\end{equation}
where 
\[
\sigma_0={\renewcommand{\arraystretch}{1}
\begin{pmatrix}1 & 0\\ 0 & 1\end{pmatrix},\;\;
\sigma_1=\begin{pmatrix}0 & 1 \\ -1 & 0 \end{pmatrix},\;\;
\sigma_2=\begin{pmatrix}0 & i \\ i & 0 \end{pmatrix},\;\;
\sigma_3=\begin{pmatrix}-i & 0 \\ 0 & i \end{pmatrix}}
\]
are matrix
representations of the units of quaternion algebra $\R_{0,2}=\BH$ :
$\e_{i}\longrightarrow\sigma_i\;(i=0,1,2),\;\e_{21}\longrightarrow\sigma_3$.
It is easy to see
that the left part of the expression (\ref{e16}) is a biquaternion
$\C_2=\C\otimes\R_{0,2}$. Recalling that $\C_2=\R_{3,0}$ and a volume element
$\omega=\e_{123}\in\R_{3,0}$ belongs to a center $\Z_{3,0}=\R\oplus i\R$,
we can write a biquaternion $X^4\e_0+X^1\e_1+X^2\e_2+X^3\e_3$, where
$\e^2_1=\e^2_2=1,\;\e_3=\e_{21}=\e_2\e_1$, in the form
\begin{multline}\label{e17}
\re X^4\e_0+\re X^1\e_1+\re X^2\e_2+\re X^3\e_3+\\
+\im X^3\e_{12}+\im X^2\e_{31}+\im X^1\e_{23}+\im X^4\e_{123}=\\
=\left(\re X^4+\omega\im X^4\right)\e_0+\left(\re X^1+\omega\im X^1\right)\e_1
+\\
+\left(\re X^2+\omega\im X^2\right)\e_2+\left(\re X^3+\omega\im X^3\right)\e_3=
\\
=X^4\e_0+X^1\e_1+X^2\e_2+X^3\e_3.
\end{multline}
Further, by means of isomorphisms $\R_{3,0}\cong\R^+_{1,3}$ and
$\R^{++}_{4,1}\cong\R^+_{1,3}\cong\R_{3,0}$ the biquaternion (\ref{e17})
may be rewritten as (like (\ref{e9})):
\begin{multline}\label{e17'}
\phi=\re X^4I+\re X^1\Gamma_{01}+\re X^2\Gamma_{02}+\re X^3\Gamma_{03}+\\
+\im X^3\Gamma_{12}+\im X^2\Gamma_{31}+\im X^1\Gamma_{23}+\im X^4\Gamma_{0123}.
\end{multline}
Or in the form (\ref{e10}) if suppose
\begin{equation}\label{e18}
\begin{array}{ccc}
\phi_1&=&\re X^4-i\im X^3,\\
\phi_2&=&\im X^2-i\im X^1,\\
\phi_3&=&\re X^3-i\im X^4,\\
\phi_4&=&\re X^1+i\re X^2.
\end{array}
\end{equation}
The formulae (\ref{e18}) define a relation between Weierstrass-Konopelchenko
coordinates and Dirac-Hestenes spinors. This relation is a direct
consequence of an isomorphism $\C_2=\C\otimes\R_{0,2}=\R_{3,0}\cong\R^+_{1,3}$.
Further, using the idempotent $\frac{1}{2}(1+\Gamma_0)\frac{1}{2}(1+
i\Gamma_{12})$ it is easy to establish (by means of (\ref{e8})) a relation
with the Dirac spinor treated as a minimal left ideal of algebra
$\R_{4,1}=\C_4\cong\M_4(\C)$:
$$\Phi=
{\renewcommand{\arraystretch}{1}
\begin{pmatrix}
\phi_1 & 0 & 0 & 0 \\
\phi_2 & 0 & 0 & 0 \\
\phi_3 & 0 & 0 & 0 \\
\phi_4 & 0 & 0 & 0
\end{pmatrix}}$$

It is obvious that we cannot directly to identify the spinor defined by
the formulae (\ref{e18}) with a generic ``physical'' spinor of electron
theory, because in accordance with (\ref{e18}) and (\ref{e16'})-(\ref{e17''})
the spinor $\phi$ depends only on two variables $z,\ov{z}$, or $x^1,x^2$
if suppose $z=x^1+ix^2$, whilst a physical spinor with four components
depends on four variables $x^1,x^2,x^3,x^4$. By this reason we will call
the spinor defined by the identities (\ref{e18}) as {\it a surface spinor},
and respectively the field $\Phi=(\phi_1,\phi_2,\phi_3,\phi_4)^T$ will be
called {\it a Dirac spinor field on surface}. The relationship between
a surface spinor $\phi(z,\ov{z})$ and a physical (space) spinor 
$\phi(x^1,x^2,x^3,x^4)$, and also a relation with the spinor representations
of surfaces in spaces $\R^{p,q}$ will be considered in a separate paper.

Further, according to Section 2, an algebraic Dirac spinor for $\R^{1,3}$
is an element of $\left\{I_{\Sigma}\right\}/\R$. Then if $\Phi_{\Sigma_0}\in
I_{\Sigma_0},\,\Phi_{\overset{.}{\Sigma}}\in I_{\overset{.}{\Sigma}}$, then
$\Phi_{\overset{.}{\Sigma}}\simeq\Phi_{\Sigma_0}\pmod{\R}$ if and only if
\begin{equation}\label{e19}
\Phi_{\overset{.}{\Sigma}}=\Phi_{\Sigma_0}u^{-1},\quad u\in\Spin_+(1,3).
\end{equation}
Here in accordance with (\ref{e8})
$$\Phi_{\Sigma_0}=\phi_{\Sigma_0}\frac{1}{2}(1+\Gamma_0)\frac{1}{2}
(1+i\Gamma_{12}).$$
The formula (\ref{e19}) defines a transformation law of the Dirac spinor.
It is obvious that a transformation group of the biquaternion (\ref{e16})
is also isomorphic to $\Spin_+(1,3)$, since
$$\Spin_+(1,3)\cong\left\{{\renewcommand{\arraystretch}{1}
\begin{pmatrix}a & c \\ b & d\end{pmatrix}\in\C_2
\,:\;\det\begin{pmatrix}a & c \\ b & d\end{pmatrix}}=1\right\}=SL(2;\C),$$
where $SL(2;\C)$ is a double covering of the own Lorentz group $\pounds_+
^{\uparrow}$.
Therefore, {\it the transformations of Weierstrass-Konopelchenko coordinates
for surfaces immersed into $\C^4$ are induced (via the relations
(\ref{e18})) transformation of a Dirac field $\Phi=(\phi_1,\phi_2,\phi_3,
\phi_4)^T$ in $\R^{1,3}$, where $\Phi\in\M_4(\C)e_{41}$ is the minimal left
ideal of $\R_{4,1}\cong\M_4(\C)$ defined in some orthonormal basis
$\Sigma\in\mathfrak{B}_{\Sigma}$.}

On the other hand, if suppose 
(following \cite{2,3,15}) that the functions $p,\psi_{\alpha}$
and $\varphi_{\alpha}$ in (\ref{e17''}) depend on the 
deformation parameter $t$, then the
deformations of $\psi_{\alpha}$ and $\varphi_{\alpha}$ are defined by a
following system:
\begin{equation}\label{e20}
\begin{array}{ccc}
\psi_{\alpha t}&=&A\psi_{\alpha}+B\varphi_{\alpha},\\
\varphi_{\alpha t}&=&C\psi_{\alpha}+D\varphi_{\alpha},
\end{array}\quad
\alpha=1,2
\end{equation}
where $A,B,C,D$ are differential operators. The equations (\ref{e20})
define integrable deformations of surfaces immersed in $\C^4$.
Let $p$ be a real-valued function;
the compatibility condition of (\ref{e20}) with (\ref{e17''}) is equivalent
to the nonlinear partial differential equation for $p$. In the simplest
nontrivial case ($A,B,C,D$ are first order operators) it is a modified
Veselov-Novikov equation \cite{31}:
$$p_t+p_{zzz}+p_{\bar{z}\bar{z}\bar{z}}+3p_{z}\omega+3p_{\bar{z}}\bar{\omega}+
\frac{3}{2}p\bar{\omega}_{\bar{z}}+\frac{3}{2}p\omega_{z}=0,$$
$$\omega_{\bar{z}}=(p^2)_z.$$
Varying operators $A,B,C,D$ one gets an infinite hierarchy of integrable
equations for $p$ (modified Veselov-Novikov hierarchy
\cite{32,33,31,2}). It is obvious that
the deformation of $\psi_{\alpha},\varphi_{\alpha}$ via (\ref{e20})
induced the deformations of the coordinates $X^i(z,\ov{z},t)$ in $\C^4$.
Moreover, according to (\ref{e16}) and (\ref{e17'}) treated as a matrix
representation of the Dirac-Hestenes spinor field $\phi$, we may say that
the mVN-deformation generates a deformation of the Dirac field
$\Phi=(\phi_1,\phi_2,\phi_3,\phi_4)^T$.
\section*{Acknowledgements}
I am grateful to Prof. B.G. Konopelchenko whose articles and preprints
are the ground for this work.
 

\begin{thebibliography}{99}
\bibitem{1} K. Weierstrass, {\it Untersuchungen \"{u}ber die Fl\"{a}chen,
deren mittlere Kr\"{u}mmung \"{u}berall gleich Null ist}, Monatsber.
Akad. Wiss. Berlin, 1866, S. 612-625.
\bibitem{2} B.G. Konopelchenko, {\it Induced surfaces and their integrable
dynamics}, Stud. Appl. Math., {\bf 96}, 9-51, (1996); preprint Institute of
Nuclear Physics, N 93-144, Novosibirsk, (1993).
\bibitem{3} B.G. Konopelchenko, {\it Multidimensional integrable systems
and dynamics of surfaces in space}, preprint of Institute of Mathematics,
Taipei, August, 1993.
\bibitem{4} R. Carroll and B.G. Konopelchenko, {\it Generalized Weierstrass-
Enneper inducing, conformal immersion and gravity}, Int. J. Modern Physics
{\bf A11}, 1183-1216, (1996).
\bibitem{5} B.G. Konopelchenko and I.A. Taimanov, {\it Generalized
Weierstrass formulae, soliton equations and Willmore surfaces}, preprint
N.187, Univ. Bochum, (1995).
\bibitem{6} B.G. Konopelchenko and I.A. Taimanov, {\it Constant mean
curvatire surfaces via an integrable dynamical system}, J. Phys. A:
Math. Gen., {\bf 29}, 1261-1265, (1996).
\bibitem{7} I.A. Taimanov, {\it Modified Novikov-Veselov equation and
differential geometry of surfaces}, Trans. Amer. Math. Soc., Ser.2,
{\bf 179}, 133-159, (1997).
\bibitem{8} I.A. Taimanov, {\it Surfaces of revolution in terms of solitons},
Annals of Global Analysis and Geometry, {\bf 15}, N 5, 419-435, (1997).
\bibitem{9} I.A. Taimanov, {\it The Weierstrass representation of closed
surfaces in $R^3$}, preprint dg-ga/9710020, 1997.
\bibitem{10} I.A. Taimanov, {\it Global Weierstrass representation and its
spectrum}, Uspechi Mat. Nauk, {\bf 52}, N 6, 187-188, (1997).  
\bibitem{11} P.G. Grinevich and M.V. Schmidt, {\it Conformal invariant
functionals of immersion of tori into $R^3$}, J. Geometry and Physics
{\bf 26},1-2, 51-78, (1998).
\bibitem{12} B.G. Konopelchenko and G. Landolfi, {\it On classical string
configurations}, Mod. Phys. Lett. A, {\bf 12}, (1997), 3161-3168.
\bibitem{13} R. Kusner and N. Schmitt, {\it The spinor representation of
surfaces in space}, preprint dg-ga/9610005, 1996.
\bibitem{14} T. Friedrich, {\it On the spinor representation of surfaces
in euclidean 3-spaces}, preprint SFB 288, N 295, TU-Berlin, 1997;
J. Geometry and Physics (to appear).
\bibitem{15} B.G. Konopelchenko and G. Landolfi, {\it Generalized
Weierstrass representation for surfaces in multidimensional Riemann spaces},
preprint math.DG/9804144, 1998; J. Geometry and Physics (to appear).
\bibitem{16} D. Hestenes, {\it Real spinor fields}, J. Math. Phys. {\bf 8},
798-808, (1967).
\bibitem{17} D. Hestenes, {\it Observables, operators, and complex numbers
in the Dirac theory}, J. Math. Phys. {\bf 16}, 556-571, (1976).
\bibitem{18} D. Hestenes, {\it Space-Time Algebra}, Gordon and Breach,
New York, 1987.
\bibitem{19} D. Hestenes and G. Sobczyk, {\it Clifford algebra to
Geometrical Calculus}, D. Heidel, Publ. Co. Dordrecht, 1984.
\bibitem{20} P. Lounesto, {\it Scalar product of spinors and an extension
of Brauer-Wall groups}, Found. Phys. {\bf 11}, 721-740, (1981).
\bibitem{21} P. Lounesto, {\it Clifford algebras and Hestenes spinors},
Found. Phys. {\bf 23}, 1203-1237, (1993).
\bibitem{22} V.L. Figueiredo, W.A. Rodrigues, Jr. and E.C. Oliveira,
{\it Covariant, algebraic, and operator spinors}, Int. J. Theor. Phys. {\bf 29},
371-395, (1990).
\bibitem{23} V.L. Figueiredo, W.A. Rodrigues, Jr. and E.C. Oliveira, {\it
Clifford algebras and the hidden geometrical nature of spinors}, Algebras,
Groups and Geometries, {\bf 7}, 153-198, (1990).
\bibitem{24} A. Crumeyrolle, {\it The primitive idempotents of the Clifford
algebras and the amorphic spinor fibre bundles}, Reports on Math. Phys. 
{\bf 25}, 305-328, (1987).
\bibitem{25} A. Crumeyrolle, {\it Orthogonal and Symplectic Clifford
Algebras}, Kluwer Acad. Publ., Dordrecht, 1991.
\bibitem{26} W.A. Rodrigues, Jr., Q.A.G. De Souza, J. Vaz, Jr., and
P. Lounesto, {\it Dirac-Hestenes spinor fields in Riemann-Cartan
spacetime}, Int. J. Theor. Phys. {\bf 35}, 1849-1900, (1996).
\bibitem{27} J. Vaz, Jr., and W.A. Rodrigues, Jr., {Maxwell and Dirac
theories as an already unified theory}, preprint hep-th/9511181, to appear
in the proceedings of the International Conference on the Theory of the
Electron, J. Keller and Oziewicz (eds.), UNAM, Mexico, (1995).
\bibitem{28} M.F. Atiyah, R. Bott and A. Shapiro, {\it Clifford modules},
Topology, {\bf 3}, (Suppl.1), 3-38, (1964).
\bibitem{29} I. Porteous, {\it Topological Geometry}, van Nostrand, London,
1969.
\bibitem{30} M. Karoubi, {\it K-Theory}, Spinger-Verlag, Berlin, 1979.
\bibitem{31} L.V. Bogdanov, {\it Veselov-Novikov equation as a natural
two-dimensional generalization of the Korteweg-de Vries equation}, Theor.
Math. Phys. {\bf 70}, 309-314, (1987).
\bibitem{32} A.P. Veselov and S.P. Novikov, {\it Finite-zone, two-dimensional
potential Schr\"{o}dinger operators. Explicit formulas and evolution
equations}, Soviet Math. Dokl. {\bf 30}, 588-591, (1984).
\bibitem{33} A.P. Veselov and S.P. Novikov, {\it Finite-zone, two-dimensional
Schr\"{o}dinger operators. Potential operators}, Soviet Math. Dokl. {\bf 30},
705-708, (1984).
\end{thebibliography}
\end{document}